\documentclass{iopart}
\usepackage{amsfonts,amssymb}        
\usepackage[dvips]{graphicx}
\usepackage[hang]{subfigure}

\def\Bbb{\mathbb}
\def\RR#1{{\Bbb R}^{#1}}
\def\RNUM{{\Bbb R}}

\def\CNUM{{\Bbb C}}
\def\ZNUM{{\Bbb Z}}
\def\NNUM{{\Bbb N}}

\def\d{\,\rmd}
\def\dx{\,\rmd x}

\def\GRAD{\nabla}
\def\LAP{\Delta}
\def\DIV{\mathrm{div}\,}

\def\PD#1#2{\frac{\partial{#1}}{\partial{#2}}}

\def\COMMA{\,,}
\def\PERIOD{\,.}
\def\SEP{{\,|\,}}
\def\VIZ#1{(\ref{#1})}
\def\NUM#1{{\tt{#1}}}

\def\REAL{\mathrm{Re}\,}

\def\BIGO{\mathrm{O}}

\def\PROOF{{\sc Proof:}\ }
\def\ATAN{\mathrm{arctan}\,}

\def\TT{\mathbb{T}}
\def\ZZ{\mathbb{Z}}

\def\KAPPA{\varkappa}

\newtheorem{theorem}{Theorem}[section]
\newtheorem{lemma}{Lemma}[section]

\newtheorem{conjecture}{Conjecture}[section]
\newtheorem{remark}{Remark}[section]
\newtheorem{definition}{Definition}[section]
\eqnobysec

\begin{document}
%
%
%
%
%
%
\title[Singular solutions]{%
Singular and regular solutions of a non-linear parabolic system
}%
\author{%
Petr Plech\'a\v{c}\dag,
Vladim\'\i{}r \v{S}ver\'ak\ddag}

\address{\dag Mathematics Institute, University of Warwick,
Coventry, CV4 7AL, UK}
\address{\ddag School of Mathematical Sciences, 
University of Minnesota, Minneapolis, MN 55455, USA}

\ead{plechac@maths.warwick.ac.uk}
\ead{sverak@math.umn.edu}

\begin{abstract}
We study a dissipative nonlinear equation modelling certain features of
the Navier-Stokes equations. We prove that the evolution of radially symmetric
compactly supported initial data does not lead to singularities in dimensions
$n\leq 4$. For dimensions $n>4$ we present strong numerical evidence supporting
existence of blow-up solutions. Moreover, using the same techniques
we numerically confirm a conjecture of Lepin regarding existence of 
self-similar singular solutions to a semi-linear heat equation.  
\end{abstract}

\submitto{Nonlinearity}


\maketitle

%
%

\maketitle
%
%
\section{Introduction}
In this paper we study solutions of the following model equation for
the time-dependent vector field $u(x,t) = (u_1(x,t),\dots, u_n(x,t))$ on
$\RR{n}\times (0,T)$
\begin{equation}\label{eq1}
   \PD{u}{t} + a\, u \GRAD u + \frac{1}{2} (1-a) \GRAD |u|^2 + 
   \frac{1}{2} (\DIV u) u = \LAP u + \KAPPA \GRAD \DIV u\COMMA
\end{equation}
where $a\in(0,1)$ and $\KAPPA\geq 0$ are given parameters.
The equation \VIZ{eq1} is of interest for various reasons. For example,
it has the same scaling properties and the same energy estimate
as the Navier-Stokes equation (NSE):  If $u(x,t)$ is a solution of \VIZ{eq1}
then also $\lambda u(\lambda x,\lambda^2 t)$ is a solution for $\lambda >0$
and, for sufficiently regular solutions with a suitable decay
at infinity, we have
%
%
\begin{eqnarray}
   \int_{\RR{n}} \frac{1}{2}|u(x,t)|^2 \dx = & 
           \int_{\RR{n}} |u(x,t')|^2 \dx \\ \nonumber
     & +\int_t^{t'}\int_{\RR{n}}
    \left[ |\GRAD u(x,s)|^2 + \KAPPA (\DIV u(x,s))^2\right]\dx \d s\PERIOD
     \label{eq2}
\end{eqnarray}
Heuristically, solutions of \VIZ{eq1} should converge to the solutions
of the NSE as $\KAPPA\to\infty$. Similar penalization schemes have 
been used in numerical algorithms for solution of 
NSE, see, e.g., \cite{Gresho}.

In  dimension $n=2$  equation \VIZ{eq1} is
``critical'' (i.e. the controlled quantities are invariant under the
scaling symmetries of the equation) and hence it is natural to expect 
that the full regularity of solutions with finite energy can be
proved by standard methods.

In this paper we shall concentrate on the super-critical case $n\geq 3$.
It is natural to expect that the theory of Leray's weak solutions applies
in this case. Moreover, it is likely that for $n=3$ 
the partial regularity results
in the spirit of Scheffer \cite{Scheffer} and Caffarelli-Kohn-Nirenberg 
\cite{Caffarelli-Kohn-Nirenberg} can be proved here as well.
We note that for $a=1/2$ the non-linear part in \VIZ{eq1} can be written
in the divergence form and consequently one can directly apply the known 
regularity theory for the NSE  in that case.
However, most questions regarding full regularity of solutions 
to \VIZ{eq1} in the case $n\ge 3$ appear to be open. 

Our aim here is to investigate the problem of finite-time blow-up
for a special class of solutions to \VIZ{eq1}.
We study solutions given by
\begin{equation}\label{eq3}
  u(x,t) = - v(r,t) x\COMMA
\end{equation}
where $r=|x|$, and $v(r,t)$ is a scalar function. Such vector fields,
usually called {\it radial vector fields}, are not divergence free unless
$v\equiv 0$ and hence the relevance of such solutions for the theory
of the NSE may be limited. Nevertheless, the behaviour of these solutions
provides an interesting insight into
various scenarios of singularity formation.

Using the radial vector field {ansatz} and substituting in the equation
\VIZ{eq1} we obtain
\begin{equation}\label{eq4}
  v_{t} = (1+\KAPPA) \left(v_{rr} + \frac{n+1}{r}v_r\right) + 3 r v v_r +
          (n+2)v^2 \COMMA
\end{equation}
where subscripts denote corresponding partial derivatives. Replacing
$v(r,t)$ by $(1+\KAPPA)v(r,(1+\KAPPA)t)$ we see that, when studying
the radial solutions, one can assume
$\KAPPA=0$ without loss of generality.

\medskip

Our first result is that in dimension $n\leq 4$ the solutions to
\VIZ{eq4} do not exhibit blow-up if there exists $C>0$
such that the initial condition 
$v(r,0)=v_0(r)$ satisfies
\begin{eqnarray} 
      -C \,\,\leq \,\,v_0(r) &\leq& C(1+r)^{-(n+2)/3}
          \;\;\mbox{when $n < 4$, and } \label{eq6}\\
      -C \,\,\leq \,\,v_0(r) &\leq& 
                             \frac{1}{r^2}\left(\frac{4}{3}\log r +C\right) 
           \;\;\mbox{when $n = 4$}\PERIOD \label{eq7}
\end{eqnarray}
On the other hand, when $v(r,0)=v_0(r) = c>0$
(where $c$ is a constant), then $v(r,t)=v(t)$ solves
\[
  \frac{dv}{dt} = (n+2) v^2\COMMA\;\;\;v(0)=c \COMMA
\]
and the solution blows up at time $t=1/(c(n+2))$.
Therefore some control
of $v_0(r)$ at infinity is necessary to prevent  formation of singularities. 

The proof, that conditions \VIZ{eq6}-\VIZ{eq7} are sufficient for preventing
blow-up, is based on an analysis of steady-state solutions to the equation \VIZ{eq4}.
The steady states can be analyzed  more or less completely since the equation
\begin{equation}\label{eq8}
  v'' + \frac{n+1}{r} v' + 3r v v' + (n+2) v^2 = 0\COMMA
\end{equation}
can be transformed to an autonomous two-dimensional dynamical system.
We briefly outline the behaviour of the steady-state solutions:
Equation \VIZ{eq8} has a solution $V\colon [0,\infty) \to (0,\infty)$
with  $V(0)=1$, $V'(0)=0$ and the following asymptotics at infinity
\begin{eqnarray*}
  V(r) &\sim & r^{-(n+2)/3}\COMMA\;\;\;\mbox{when $1 < n < 4$} \\
  V(r) &\sim & r^{-2}\left(\frac{4}{3}\log r + C\right)\COMMA\;\;\;\mbox{when $n= 4$}\\
  V(r) &\sim & r^{-2} \COMMA\;\;\;\mbox{when $n>4$}\PERIOD
\end{eqnarray*}
Using the scaling symmetry we obtain a one-parameter family of solutions
\[
  v_\lambda(r) = \lambda^2 V(\lambda r)\PERIOD
\]

There are also other interesting steady-state solutions. It turns out
that even for radial solutions,  {\it weak} solutions of
\VIZ{eq1} can exhibit the following non-trivial behaviour:
\begin{description}
  \item[{\rm (i)}] formation of singularities with a different 
                   rate of blow-up than suggested by scaling,
  \item[{\rm (ii)}] violation of local energy inequality,
  \item[{\rm(iii)}] significant non-uniqueness.
\end{description}
We present more specific discussion of these phenomena in Section~\ref{sec2}.

The second group of the results, we shall discuss, concerns the blow-up
behaviour of solutions to the equation \VIZ{eq4} in dimensions $n>4$.
In this case our results are based on combination of
analytical arguments and numerical calculations.
We will present  strong
evidence that for $n>4$ and suitable compactly supported initial data
there exist solutions of \VIZ{eq4} that
form a singularity in finite time.

A natural class of singularities for the equation \VIZ{eq4} are
 self-similar singularities of the type
\[
  v(r,t) = \frac{1}{2\kappa (T-t)} 
           w\left(\frac{r}{\sqrt{2\kappa(T-t)}}\right)\COMMA
\]
where $\kappa>0$, $T>0$ are parameters and $w$ is a function defined
on $[0,\infty)$. The equation for $w$ is
\begin{equation}\label{eq9}
  w'' + \frac{n+1}{r} w' -\kappa r w' + 3 r w w' + (n+2) w^2 - \kappa w =0\COMMA
\end{equation}
together with the natural boundary conditions
\begin{eqnarray}
  w(0) &=& \alpha >0\COMMA\;\;\;w'(0)=0\COMMA  \label{eq10a} \\
  w(r) &\sim & r^{-2}\COMMA\;\;\mbox{as $r\to\infty$.} \label{eq10b}
\end{eqnarray}
The scaling symmetry $(w,\kappa) \to (\lambda^2 w(\lambda r),\lambda^2 \kappa)$
allows us to fix $\kappa$ and vary only $\alpha$. There are no non-trivial
solutions of \VIZ{eq9}-\VIZ{eq10b} for $n\leq 4$. We conjecture that for
$n>4$ and a fixed $\kappa$, the number of non-trivial solutions is determined
by the number of roots of the solution to a linearization
of \VIZ{eq9} around its trivial equilibrium
\[
  \bar w = \frac{2\kappa}{n+2}\PERIOD
\]
The relevant linear problem is 
\begin{eqnarray}
   && z'' + \left( \frac{n+1}{r} - \kappa r \frac{n-4}{n+2}\right) z'
       + 2\kappa z = 0  \label{eq11}\\
   && z(0)=1\COMMA\;\; z'(0)=0\PERIOD \nonumber
\end{eqnarray}
The solution of this problem can be written  explicitly in terms of 
a confluent hyper-geometric function:
\[
   z(r) = M\left(-\frac{n+2}{n-4},\frac{n+2}{2},\kappa\frac{r^2}{2}\right)
                                                   \PERIOD
\]
We refer the reader to \cite{SpecFun} for basic properties of hypergeometric functions.
It turns out that $z(r)$ is a polynomial whenever $(n+2)/(n-4)$ is an integer.
The number of zeros $m=m(n)$ of $z(r)$ in $(0,\infty)$ is
the smallest integer greater or equal to $(n+2)/(n-4)$. 
We recall that we assume $n>4$ at this point.
We conjecture that, for a fixed $\kappa$ and $n>4$, the number of 
solutions of \VIZ{eq9}-\VIZ{eq10b} is $m(n)-2$. 

There are certain  similarities between the behaviour of solutions 
of \VIZ{eq4} and solutions of a system arising in mathematical biology studied
in \cite{Brenner}.

\medskip

In the last section we study some related problems concerning a widely
studied  semi-linear
heat equation. There appear to be striking similarities between the
formation of singularities in the equation \VIZ{eq4} for $n>4$ and 
similar behaviour of non-negative radial solutions of 
\begin{equation}\label{eq9a}
   \PD{v}{t} = \LAP v + v^{2\sigma + 1}\COMMA\;\;\;\mbox{in $\RR{n}\times(0,t)$}
\end{equation}
for $n>10$ and $\sigma > \sigma_c(n)=2/(n-4-2\sqrt{n-1})$. For the significance
of the critical exponent $\sigma_c(n)$ see, for example, 
\cite{G-Ni-Wang,He-Va}.

Self-similar singular solutions of \VIZ{eq9a} are given by
\[
  v(x,t) = \left(2\kappa(T-t)\right)^{-1/\sigma} 
          w \left(\frac{|x|}{\sqrt{2\kappa(T-t)}}\right)\COMMA
\]
where $\kappa,T >0$ are parameters and $w$ is a function on $[0,\infty)$.
The function $w$ is a solution of the boundary value problem
\begin{eqnarray}
   && w'' + \frac{n-1}{r}w' - \kappa \left(\frac{w}{\sigma} + r w'\right)
         + w^{2\sigma+1} = 0 \label{eq13} \\
   && w(0)=\alpha\COMMA\;\;\;w'(0)=0\COMMA\;
       \mbox{and $w(r)\sim r^{-1/\sigma}$ as $r\to\infty$} \label{eq14}
\end{eqnarray}
The existence of non-trivial solutions to \VIZ{eq13}-\VIZ{eq14} depends
on $n$ and $\sigma$ in the following way. If
\begin{eqnarray*}
  \fl \;\;\;\;\; \sigma \leq \frac{2}{n-2}\COMMA\;
             \mbox{there are no non-trivial solutions, see \cite{Kohn-Giga},} \\
  \fl \;\;\;\;\; \frac{2}{n-2} < \sigma < \frac{2}{n-4-2\sqrt{n-1}}\COMMA\;
             \mbox{there are infinitely many solutions, see \cite{Troy},} \\
  \fl \;\;\;\;\; \frac{2}{n-2} < \sigma < \frac{3}{n-10}\COMMA\;
             \mbox{there exists at least one non-trivial solution, see \cite{Lepin}.}
\end{eqnarray*}
One of the open problems for \VIZ{eq13}-\VIZ{eq14} is to determine 
the exact range
of parameters for which the boundary-value problem has a non-trivial
solution.
In Section~\ref{heq} we present strong evidence that the sufficient
condition of Lepin (\cite{Lepin}) $\sigma < 3/(n-10)$ is also necessary
for the existence of non-trivial positive solutions of \VIZ{eq13}-\VIZ{eq14}.

\section{Solutions of the model equation}\label{sec2}

\subsection{Phase portrait}\label{subsec1}
In this section we analyze the steady-state solutions of the equation
\begin{equation}\label{SS}
  v'' + \frac{n-1}{r}v' + 3 r v v' +(n+2) v^2 = 0 \PERIOD
\end{equation}
The invariance of solutions to \VIZ{SS}
under the scaling $v(r) \to \lambda ^2 v(\lambda r)$
suggests the change of variables:
\[
  v = r^{-2} w \COMMA\;\;\; r = e^s\COMMA
\]
which transforms \VIZ{SS} into an autonomous equation
\begin{equation}\label{AE}
  w'' + 3w w' - (4-n)w^2 + 2(n-2)w - (4-n)w' = 0\PERIOD
\end{equation}
With a slight abuse of notation the prime $'$ now denotes differentiation 
with respect to the new independent variable $s$.

\begin{figure}[ht]
\caption{Heteroclinic orbits in different dimensions $n$.}\label{fig1}
\subfigure[\label{pgt4} $n>4$.]{
\begin{minipage}{6.5cm}
\begin{center}
\includegraphics[width=6cm]{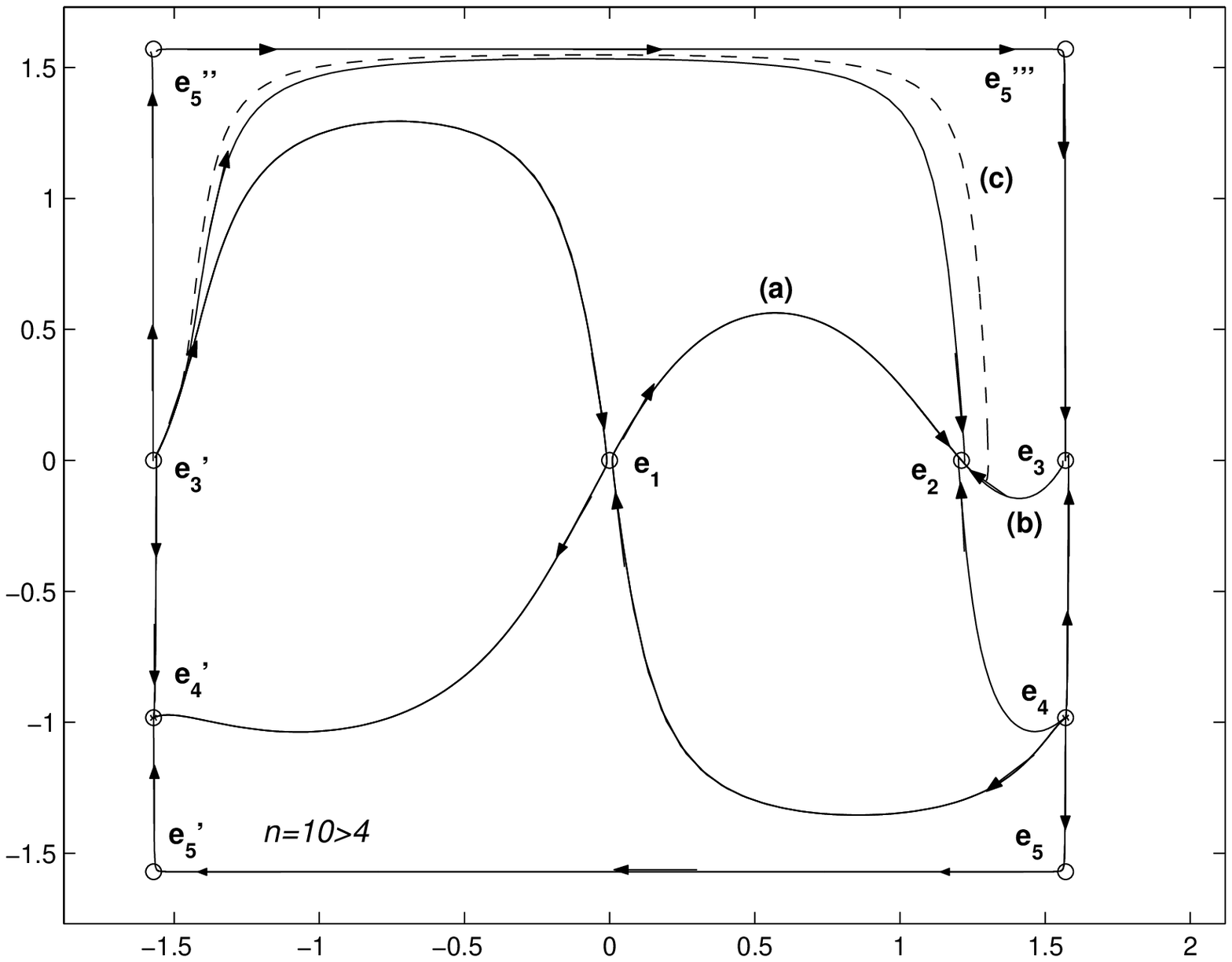}
\end{center}
\end{minipage}}
\subfigure[\label{peq4} $n=4$.]{
\begin{minipage}{6.5cm}
\begin{center}
\includegraphics[width=6cm]{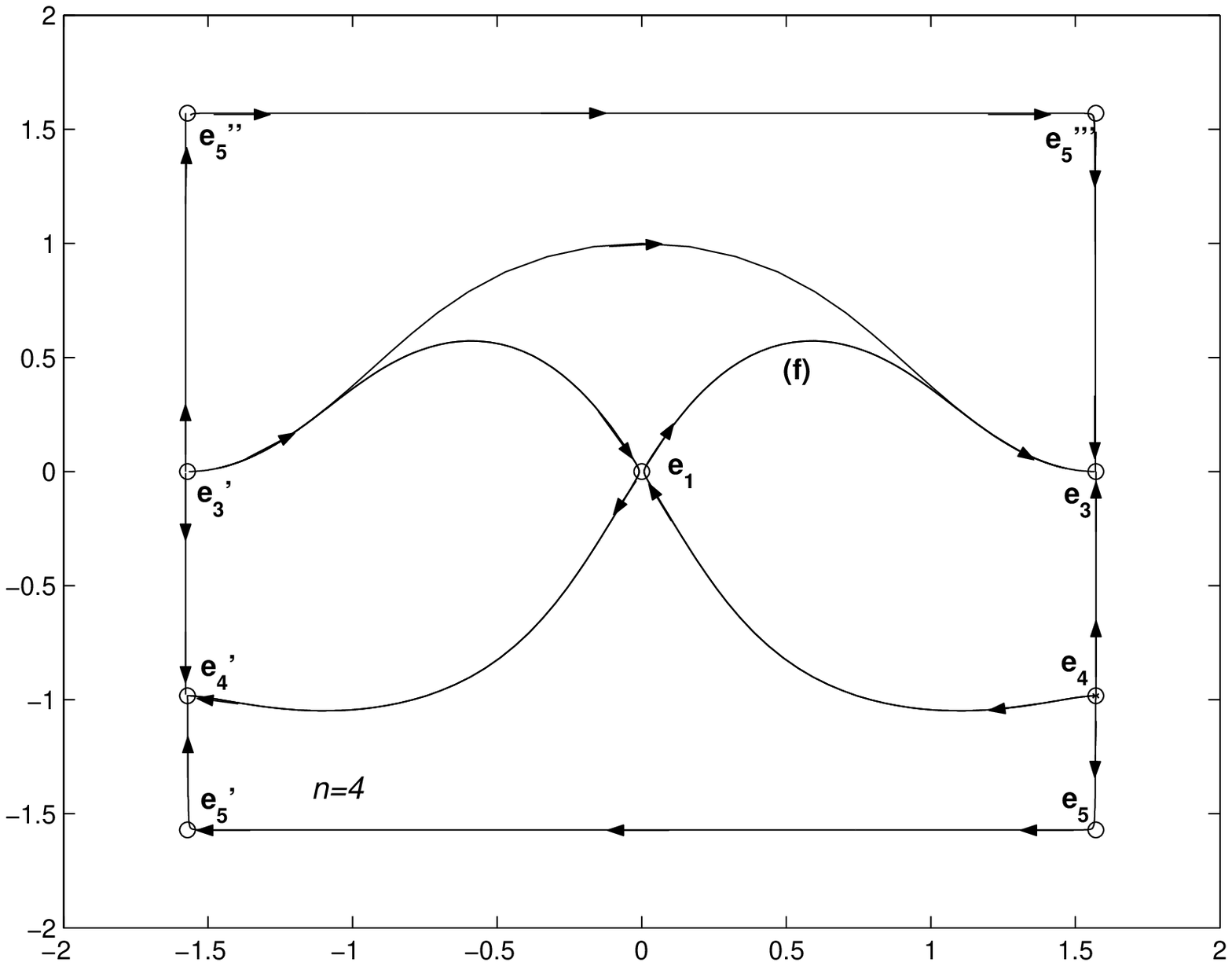}
\end{center}
\end{minipage}}\\

\subfigure[\label{plt4} $n<4$.]{
\begin{minipage}{13cm}
\begin{center}
\includegraphics[width=6cm]{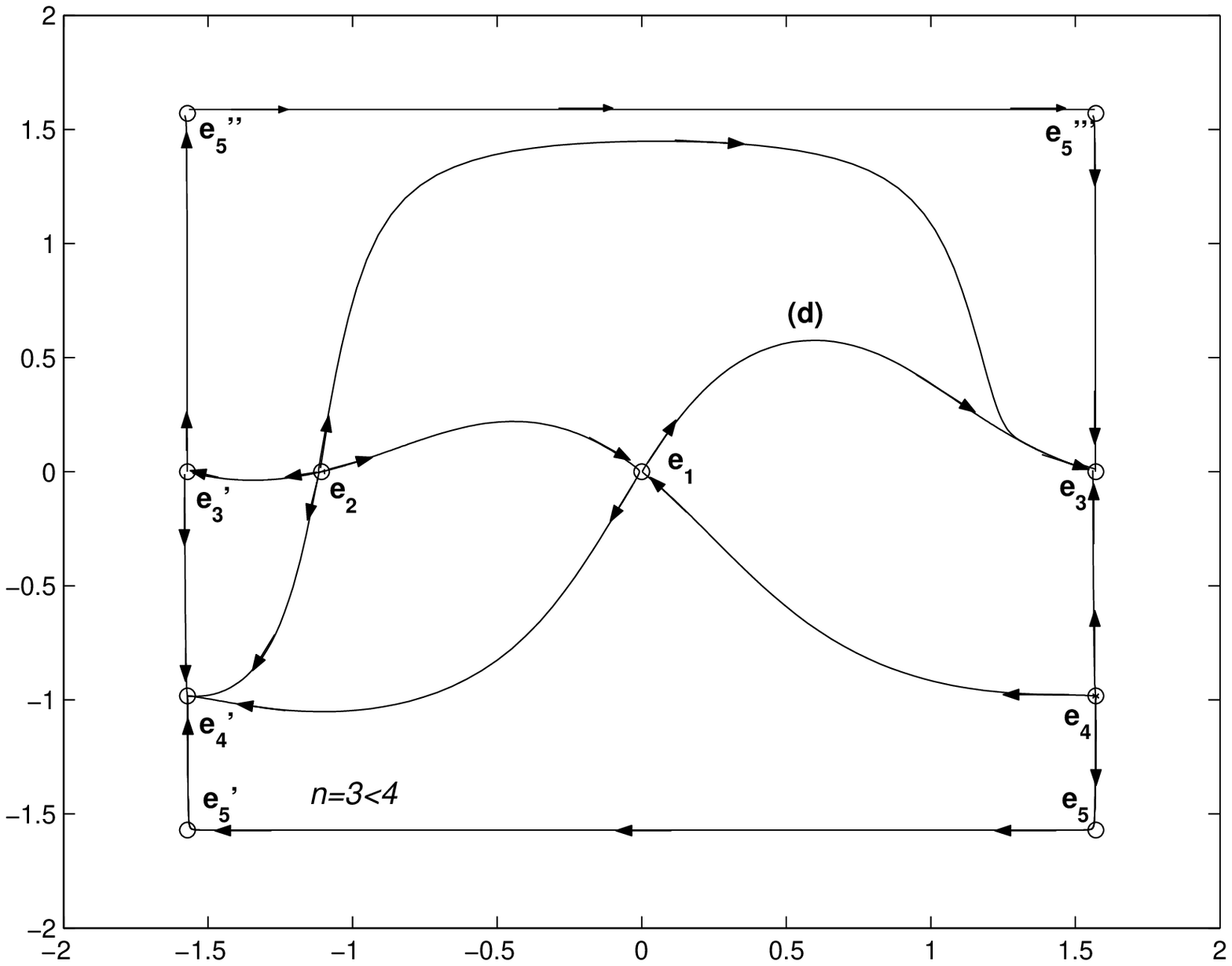}
\end{center}
\end{minipage}}
\end{figure}

We are interested in the phase portrait of the 
vector field in $\RR{2}$ defined by \VIZ{AE}. 
Global properties are  best studied in a suitable 
compactification of $\RR{2}$.
It turns out that the transformation of variables
\[
  w(s) = \tan(\phi(s))\COMMA\;\; w'(s) = \frac{\tan\psi(s)}{\cos^2\phi(s)}\COMMA
\]
leads to a compactification which works well in the case at hand. 
In the new variables the equation \VIZ{AE} becomes
\begin{eqnarray}
  \frac{d\phi}{ds} & = & \frac{1}{\cos\phi \cos\psi} P(\phi,\psi)
                         \COMMA \label{ET} \\
  \frac{d\psi}{ds} & = & \frac{1}{\cos\phi \cos\psi} Q(\phi,\psi)
                         \COMMA \nonumber
\end{eqnarray}
where 
\begin{eqnarray*}
\fl P(\phi,\psi) =  \cos\phi \sin\psi\COMMA \\
\fl Q(\phi,\psi) =  -\sin\psi\cos\psi\sin\phi (2 \sin\psi + 3\cos\psi) 
                    +2(n-2) \sin\phi\cos^2\phi\cos^3\phi \\
                 \lo+ (4-n) \cos^2\psi\cos\phi(\sin^2\phi\cos\psi + \sin\psi) 
                 \PERIOD
\end{eqnarray*}
Hence the integral curves of the vector field \VIZ{ET} on the torus 
$\TT^2 = \RR{2}/_{2\pi\ZZ^2}$ are defined by
\begin{equation}\label{MET}
  \frac{d\phi}{ds}  =  P(\phi,\psi)\COMMA\;\;\;\;\;\;
  \frac{d\psi}{ds}  =  Q(\phi,\psi)\PERIOD
\end{equation}
Due to the periodicity we have
\begin{eqnarray*} 
&& P(\phi+\pi,\psi) = P(\phi,\psi+\pi) = -P(\phi,\psi) \\
&& Q(\phi+\pi,\psi) = Q(\phi,\psi+\pi) = -Q(\phi,\psi)\PERIOD
\end{eqnarray*}
Therefore it is sufficient to analyze
the flow in the region $|\phi|\leq\pi/2$, $|\psi|\leq \pi/2$.
We proceed with a description of equilibria and important
heteroclinic orbits in the region $(-\pi/2,\pi/2]\times(-\pi/2,\pi/2]$.
Both, the orbits and the equilibria are computed in a fully rigorous way
by standard methods.
All the other equilibria are obtained  by shifts along the 
coordinate axes by $k\pi$.
From the information about the equilibria and the heteroclinic orbits
we determine the full phase portrait of the system \VIZ{MET}. 
The phase portraits differ for different values of the dimension and
we sketch the three different cases $n<4$, $n=4$ and $n>4$.
%
\begin{table}
\caption{Equilibria of the system \VIZ{MET} in the region 
         $(-\pi/2,\pi/2]\times(-\pi/2,\pi/2]$.} \label{TAB1}
\begin{indented}
\item[] \begin{tabular}{@{}ll}
\br
       {Equilibrium $(\phi,\psi)$}      &  {Linearization}          \\ 
                                        &  eigenvalues, eigenvectors \\
\mr
$e_1$: $(\phi_1,\psi_1) = (0,0)$ & $\lambda_1 = 2$, 
                                  $\left(\begin{array}{c}
                                         1 \\
                                         2
                                    \end{array}\right)\COMMA$
                                 $\lambda_2 = -(n-2)$, 
                                 $\left(\begin{array}{c}
                                         -1 \\
                                         n-2
                                  \end{array}\right)$        \\
$e_2$: $(\phi_2,\psi_2) = (\ATAN \frac{2(n-2)}{n-4},0)$ 
                                & $\lambda_{1,2} >0$ if $n<4\COMMA$
                                  $\lambda_{1,2} <0$ if $n>4$,          \\                      
$e_3$: $(\phi_3,\psi_3) = (\frac{\pi}{2},0)$ & $\lambda_1 = -3$, 
                                  $\left(\begin{array}{c}
                                         0 \\
                                         1
                                    \end{array}\right)\COMMA$
                                  $\lambda_2 = 0$, 
                                  $\left(\begin{array}{c}
                                         3 \\
                                         n-4
                                    \end{array}\right)$                   \\
$e_4$: $(\phi_4,\psi_4) = (\frac{\pi}{2},\ATAN (-\frac{2}{3}))$
                                & $\lambda_1 = 1\COMMA$
                                  $\lambda_2 = \frac{39}{4}\cos^3\psi_4 > 1$  \\
$e_5$: $(\phi_5,\psi_5) = (\frac{\pi}{2},\frac{\pi}{2})$ 
                                & $\lambda_1 = -1$,
                                   $\left(\begin{array}{c}
                                         1 \\
                                         0
                                    \end{array}\right)\COMMA$ 
                                  $\lambda_2 = 2\COMMA$
                                   $\left(\begin{array}{c}
                                         0 \\
                                         1
                                    \end{array}\right)$
   \\ 
\br        
\end{tabular}
\end{indented}
\end{table}
Important heteroclinic connections in a typical phase portrait when $n>4$ are
depicted in Figure~\ref{fig1}(a),                                    
where $n=5$ was used in numerical computations.
Similarly Figure~\ref{fig1}(b)                                       
depicts the phase portrait for $n=4$ and
Figure~\ref{fig1}(c)                                                 
for $n<4$ (computed for $n=3$).
The curves denoted by
$a$, $b$, $c$, $d$, $f$ are important in our considerations and
their meaning is explained below.
%
%

\subsection{Specific solutions}\label{subsec2}
Analysis in the previous section enables us to show existence of
some important solutions of \VIZ{SS}.

\medskip

\noindent
{\bf Orbit (a)} 
The heteroclinic orbit (a) in Figure~\ref{fig1}(a)                    
gives immediately the solution $v:[0,\infty) \to \RNUM$, $v>0$,
\begin{eqnarray*}
         && v(0)=1\COMMA\;v'(0)=0\COMMA \\
         && v(r) \sim r^{-2}\COMMA\;\;\mbox{as $r\to\infty$} \PERIOD
\end{eqnarray*}

\noindent
{\bf Orbit (b)} 
The orbit (b) in the phase portrait for $n>4$ (Figure~\ref{fig1}(a))   
corresponds to the solution 
$v:[0,\infty) \to \RNUM$, $v>0$,
\begin{eqnarray*}
    && v(r) \sim r^{-\frac{n+2}{3}} \COMMA\;\;\mbox{as $r\to 0_+$}\COMMA\\
    && v(r) \sim r^{-2}\COMMA\;\;\mbox{as $r\to\infty$} \PERIOD
\end{eqnarray*}
The behaviour as $r\to 0_{+}$ is determined from the slope of the center 
manifold at the equilibrium point $e_3$. If $(\phi,\psi)$ approaches $e_3$ along
(b) we have the following asymptotic expansion
\begin{equation}\label{ASEXP}
      \psi = \frac{n-4}{3}(\phi - \frac{\pi}{2}) + a_2 (\phi - \frac{\pi}{2})^2
                   + \dots + a_k  (\phi - \frac{\pi}{2})^k \dots
\end{equation}
The series is not convergent as is typical for similar situations
involving center manifolds. The coefficients $a_k$ can be calculated
but we will not need their exact values for the subsequent analysis.
The expansion \VIZ{ASEXP} gives $w(s)\to \infty$ as $s\to -\infty$ and
\begin{equation}\label{ASEXP2}
       w'(s) = -\frac{n-4}{3} w(s) + b_0 + b_{-1} w^{-1}(s) + \dots + 
                       b_{-k} w^{-k}(s) \dots\COMMA
\end{equation}
which again is an asymptotic expansion, not necessarily convergent.
From  \VIZ{ASEXP2} we can see the described behaviour of $v(r)$. 

\noindent
{\bf Orbit (c)}
A similar analysis of the orbit (c) in Figure~\ref{fig1}(a)              
gives the solution $v:[0,\infty) \to \RNUM$
\begin{eqnarray*}
      && v(r) \sim r^{-\frac{n+2}{3}} \COMMA\;\;\mbox{as $r\to 0_+$}\COMMA\\
      && v(r) \sim r^{-2}\COMMA\;\;\mbox{as $r\to\infty$} \PERIOD
\end{eqnarray*}
 
\noindent
{\bf Orbit (d)}
Exploiting again the known slope of the center manifold at $e_3$
one easily sees that the heteroclinic orbit (d) in the phase portrait
for $n<4$ corresponds to the solution $v:[0,\infty) \to \RNUM$
\begin{eqnarray*}
     && v(0)=1 \COMMA \;\;v'(0) = 0\COMMA \\
     && v(r) \sim r^{-\frac{n+2}{3}}\COMMA\;\;\mbox{as $r\to\infty$} \PERIOD
\end{eqnarray*}

\noindent
{\bf Orbit (f)}
In the case $n=4$ we have an explicit equation for the solution
corresponding to the orbit (f)
\[
    -\frac{1}{3}w' + \frac{4}{9} \log \frac{1}{1 - 3 w'/4} = \frac{w^2}{2}\COMMA
\]
which leads to the solution $v:[0,\infty) \to \RNUM$
\begin{eqnarray*}
     && v(0)=1 \COMMA \;\;v'(0) = 0\COMMA \\
     && v(r) \sim r^{-2}\left(\frac{4}{3}\log r + C\right)
         \COMMA\;\;\mbox{as $r\to\infty$} \PERIOD
\end{eqnarray*}

\medskip

The existence of solutions (d) and (f) implies the following result:
\begin{theorem}
Assume that $u_0(x)=-v_0(|x|) x/|x|$ is a smooth radial vector field on $\RR{n}$ and
assume that one of the following conditions is satisfied \\
(i)  $n=2,3$ and 
\[
  -C\,\,\le\,\, v_0(|x|) \leq C (1 + |x|)^{-(n-1)/3}\COMMA 
             \;\;\mbox{for some $C>0$,} 
\]
(ii) $n=4$ and 
\[
  -C\,\,\le\,\, v_0(|x|) \leq \frac{4/3 \log(1 + |x|) + C}{1 + |x|}
             \COMMA\;\;\mbox{for some $C>0$,} 
\]
(iii) $n\geq 5$ and 
\[
  -C\,\,\le\,\, v_0(|x|) \leq \gamma \frac{(n-2)}{(n-4)}\frac{2}{|x|}\COMMA
                   \;\;\mbox{for some $\gamma < 1$ and $C>0$.}
\]
Then the equation \VIZ{eq1} with $\KAPPA = 0$ has a global bounded solution
with the initial condition $u_0$.
\end{theorem}

\noindent
\PROOF  The proof of the above theorem follows from
the analysis of the two-dimensional system discussed above.
We use the solutions constructed in this analysis as barriers
in the equation \VIZ{eq4}. 

\begin{remark}
As we have noted in the introduction, when studying radial solutions
to the equation \VIZ{eq4}, one can consider $\KAPPA =0$ without loss of
generality.
\end{remark}

\subsection{Singular steady states for $n\ge 5$}

We consider the following boundary-value problem for radial
vector fields $u:\RR{n}\to\RR{n}$:
\begin{eqnarray}
&& - \LAP u + \frac{1}{2} u \GRAD u + 
     \frac{1}{2}\GRAD\frac{|u|^2}{2} +  \frac{1}{2}u \DIV u = 0 \COMMA
     \;\;\mbox{in $B_1$} \label{BV} \\
&& u(x) = - b\, x \COMMA\;\;\mbox{on $\partial B_1$}\COMMA \label{BVa}
\end{eqnarray}
where $B_1 = \{ x\in\RR{n}\SEP |x|\leq 1\}$ denotes the unit ball 
and $b\in \RNUM$.

\begin{definition}
  We say that $u$ is a {\rm weak solution}
  of \VIZ{BV} if it belongs to $W^{1,2}(B_1)$, satisfies the equation
  in the sense of distributions, and the boundary condition is
  satisfied in the sense of traces.

  \noindent
  Furthermore, we say that $u$ is a {\rm suitable weak solution}
  of \VIZ{BV} if it is a weak solution and satisfies a local version of the 
  energy inequality
  \begin{equation}\label{EI}
     \int_{B_1} \left[ |\GRAD u|^2\phi - \frac{1}{2} |u|^2 \LAP \phi -
         \frac{1}{2}(u\cdot\GRAD\phi)|u|^2\right] \leq 0\COMMA
  \end{equation}
  for each smooth $\phi\geq 0$ compactly supported in $B_1$.
\end{definition}

It is easy to see that radial weak solutions are smooth away from the
origin. For solving \VIZ{BV}-\VIZ{BVa} we can use orbits (a) and (c)
discussed in the previous section. 
The orbit (a) gives a smooth solution when $0\le b<2(n-2)/(n-4)$.
When $b\ge2(n-2)/(n-4)$, problem \VIZ{BV}-\VIZ{BVa} does not have 
a smooth radial solution. 
In that case a {\it suitable weak solution} as defined
above can be obtained from the orbit (c). 
It is interesting to note that this solution has a singularity at the 
origin which is not asymptotically self-similar. The self-similar rate of
blow-up would be $|x|^{-1}$, whereas the actual rate is 
$|x|^{-(n-1)/3}$. Moreover, we note that the solution satisfies
the inequality \VIZ{EI}, but it does not satisfy the 
{\it local energy identity}, i.e., \VIZ{EI} with equality.

Apart from the suitable weak solution described above
there exist many other weak solutions to \VIZ{BV} when
 $b> 2(n-2)/(n-4)$. These solutions are constructed
using heteroclinic orbits such as the orbit (c) in the phase
portrait for $n>4$. However, these solutions  {\it do not}
satisfy \VIZ{EI}.

The above discussion does not cover the case $b<0$, which is left
to the reader as an exercise. 

\subsection{Self-similar singular solutions}

In this section we study radial self-similar solutions of \VIZ{eq1}, i.e.
 solutions of the form
\begin{equation}\label{S}
   u(x,t) = -\frac{1}{2\kappa(T-t)} w\left( \frac{|x|}{\sqrt{2\kappa(T-t)}}\right) x \COMMA
\end{equation}
where $T\in \RNUM$ and $\kappa> 0$. Assuming (without loss of generality) that $\KAPPA=0$
in the equation \VIZ{eq1}, we obtain the equation
\begin{equation}\label{SE}
   w'' + \frac{n+1}{r} w' - \kappa r w' + 3 r w w' + (n+2) w^2 - \kappa w = 0 \PERIOD
\end{equation}

As is usual in similar situations, it is useful to interpret the equation \VIZ{SE}
as an equation of motion for a particle with the unit mass which is moving in a potential field,
given by the potential $V(w) = (n+2)/3 w^3 - \kappa w^2/2$, in the presence of damping
$\mu(r,w) w'$ where $\mu(r,w) = (n+1)/r + 3 r w - \kappa r$. Using this notation we can write 
\[
   w'' + \mu(r,w) w' = -\frac{\partial}{\partial w} V(w)\PERIOD
\]
As in the case of steady-state solutions it is useful to introduce
new variables
by the transformation $w(r) = r^{-2} u(r)$ and $r = e^{s}$. 
In the new variables we have the equation
\begin{equation}\label{SEE}
   u'' + (n-4) u' + (n-4) u^2 - 2(n-2) u + 3u u' = \kappa e^{2s} u' \COMMA
\end{equation}
or 
\begin{equation}\label{SEEbis}
   u'' + \tilde\mu(s,u) u' = -\frac{\partial}{\partial u} \tilde V(u)\COMMA
\end{equation}
where
\[
  \tilde \mu(s,u) = n-4 + 3 u + \kappa e^{2s}\COMMA\;\;\;
  \tilde V(u) = \frac{n-4}{3} u^3 - (n-2) u^2\PERIOD
\]
Natural boundary conditions for \VIZ{SE} in the context of self-similar
singularities are
\begin{eqnarray}
   w(0) & = & \alpha\in \RNUM\COMMA\;\;\; w'(0) = 0 \;\;\mbox{and} \label{BC1}\\
   w(r) & \sim & r^{-2}\COMMA\;\;\;\mbox{as $r\to \infty$}\PERIOD\label{BC2}
\end{eqnarray}
The  condition at infinity together with \VIZ{S} leads to the blow-up profile
\[
  u(x,T) = - c \frac{x}{|x|^2}\COMMA\;\;\mbox{for some $c\in\RNUM$}\PERIOD
\]
One can easily see that $c>0$ in this case. 

The boundary-value problem \VIZ{SE}, \VIZ{BC1}-\VIZ{BC2} is studied
as a non-linear eigenvalue problem with parameters $\alpha$ and $\kappa$
considered as unknowns. Due to the scaling symmetry $(\alpha,\kappa,w) \to 
(\lambda^2 \alpha,\lambda^2\kappa,\lambda^2 w(\lambda r))$ we can assume
that $\kappa = 1$ without loss of generality.

Our analysis and numerical computations support the following conjecture
\begin{conjecture}
   Suppose  we define
   \[
      \nu(n) = \min \left\{ k\in \NNUM\SEP k\geq \frac{n+2}{n-4}\right\}\COMMA
   \]
   then for a fixed $\kappa > 0 $ and $ n>4$ the boundary-value problem
   \VIZ{SE}, \VIZ{BC1}-\VIZ{BC2} has $\nu(n) - 2$ non-trivial solutions.
\end{conjecture}

\begin{figure}[ht] 
  \caption{\label{fig2} Plots of $u(r)=-w(r)$ where $w(r)$ solves \VIZ{eq9}-\VIZ{eq10b}
    in dimension $n=5$ and $\kappa=7/2$ (so that the trivial equilibrium is $\bar w=1$). 
    The inset shows the behaviour
    of solutions around the equilibrium $\bar w = 1$.}
  \begin{center}
     \includegraphics[width=.8\linewidth]{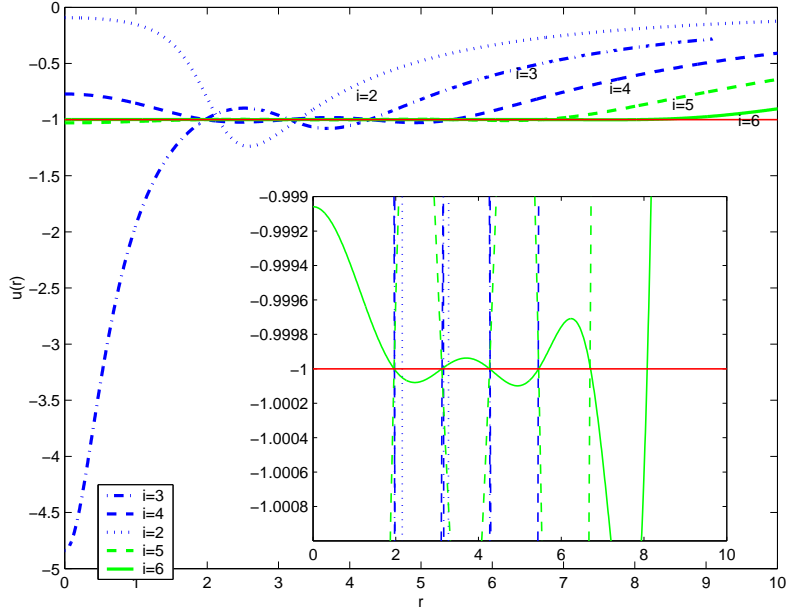}
  \end{center}
\end{figure}

\begin{remark}
  \begin{description}
      \item[{\rm (i)}] Results of Section~\ref{subsec2} imply that \VIZ{SE}, \VIZ{BC1}-\VIZ{BC2}
                                       has no non-trivial solution for $n\leq 4$.
      \item[{\rm(ii)}] Since we impose the boundary conditions \VIZ{BC1} and \VIZ{BC2}
                                       the solutions $\bar w = \frac{2\kappa}{n+2}$ and 
                                       $w(r) = 2\frac{(n-2)}{(n-4)}r^{-2}$
                                       are excluded from the count.
      \item[{\rm(iii)}] The solutions for $n=5$ are depicted in Figure~\ref{fig2}
   \end{description}
\end{remark}

The following simple observation is useful.
\begin{lemma}\label{lemma1}
  Every solution of the initial value problem \VIZ{SE} and \VIZ{BC1} defined 
  on $(0,\infty)$ is positive on $(0,\infty)$.
\end{lemma}

\noindent
\PROOF The lemma follows easily from the interpretation of the 
problem using the equation \VIZ{SEEbis}. 

\medskip

We let $w_\alpha$ be the solution of the initial-value problem \VIZ{SE}
with the initial conditions \VIZ{BC1} and $\kappa=1$. 
Furthermore we define
\[
  R_\alpha = \sup\{ R\in(0,\infty)\SEP w_\alpha\geq 0 \COMMA\mbox{ on}\, (0,R)\}\PERIOD
\]
By Lemma~\ref{lemma1} a necessary condition that $w_\alpha$ is a solution to
the boundary-value problem \VIZ{SE}, \VIZ{BC1}-\VIZ{BC2} is that $R_\alpha = +\infty$.
Moreover, it is likely that, except for perhaps some special values of $n$, 
this condition is also sufficient. The reason will become
apparent later.
For $n>4$ we let $\bar \alpha_n = 2/(n+2)$ and $\beta_n = 2(n-2)/(n-4)$ and
we define $w_\infty = \beta_n/r^2$. Note that this function satisfies \VIZ{SE}
and also $w_{\bar\alpha_n} = \bar\alpha_n$. 

The following formal calculation suggests that $w_\alpha\to w_\infty$ 
as $\alpha\to 0$. It is convenient to work
with the equation \VIZ{SEE}. In the coordinates $(s,u)$ the condition
\VIZ{BC1} becomes
\begin{equation} \label{BCE1}
   \lim_{s\to -\infty} e^{-2s} u(s) = \alpha \COMMA\;\;\mbox{and}\;\; 
   \lim_{s\to -\infty} e^{-3s} (u'(s) - 2 u(s)) = 0 \PERIOD
\end{equation}
We denote by $U(\alpha,\kappa,s)$ the solution of \VIZ{SEE} and \VIZ{BCE1}
and by $(0,s_{\alpha,\kappa})$ the maximal interval where $s\to U(\alpha,\kappa,s)$
is defined.
We note that $U(e^{2s'}\alpha,\kappa,s) = U(\alpha,e^{-2s'}\kappa,s+s')$ and will
investigate the behaviour of $U(\alpha,e^{-2s'}\kappa,s+s')$ as $s'\to\infty$. 
The linearization of \VIZ{SEE} at the equilibrium $u(s)=\beta_n$, corresponding
to $\bar w_n$ in the coordinates $(r,u)$, is defined by the equation
\begin{equation}\label{EZ}
z'' + \left( n-4 + \frac{3}{\beta_n} \right)  z' + 2(n-2)  z = \kappa e^{2s} z'\PERIOD
\end{equation}
It is easy to see that the general solution of \VIZ{EZ} is given in terms of confluent
hypergeometric functions by
\begin{equation}\label{GS}
z(s) = C_1 e^{\lambda_1 s} M(a_1, b_1, \kappa \frac{1}{2} e^{2s}) +
           C_2 e^{\lambda_2 s} M(a_2, b_2, \kappa \frac{1}{2} e^{2s})\COMMA
\end{equation}
where $a_j = \lambda_j/2$, $b_j=1 + \lambda_j + n - 4 + 3\beta_n$ and the
constants $\lambda_2 < -1 < \lambda_1 < 0$ are roots of the characteristic polynomial for the
linear second-order differential operator that defines the left hand side of \VIZ{EZ}.
The function $M$ is one of the standard confluent hyper-geometric functions, see, e.g., \cite{SpecFun}.
For small $\kappa>0$ and large $s$ such that $\kappa e^s$ is controlled,
the formula \VIZ{GS} leads to
\begin{equation}\label{APP}
U(\alpha,\kappa,s) = \beta_n - \gamma  e^{\lambda_1 s} M(a_1, b_1, \frac{1}{2}\kappa e^{2s}) + 
                                             \mbox{higher order terms}\PERIOD
\end{equation}
We use the approximation \VIZ{APP} to study asymptotic behaviour of $U(\alpha,e^{-2s'}\kappa, s+s')$
for a fixed $s$ and $s'\to \infty$. This is possible since the second term of \VIZ{APP}
becomes $-\gamma e^{\lambda_1 (s+s')} M(a_1,b_1,\frac{1}{2}\kappa e^{2s})$ which is small for large
$s'$ and fixed $s$, since $\lambda_1<0$.
Hence we conclude that
\begin{equation}\label{ASY}
  \fl U(\alpha e^{2s'},\kappa,s) = U(\alpha,e^{-2s'}\kappa, s+s') = 
                             \beta_n - e^{\lambda_1 s'} \gamma e^{\lambda_1 s} 
                             M(a_1,b_1,\frac{1}{2}\kappa e^{2s}) + \mbox{h.o.t}\PERIOD
\end{equation}
This calculation formally shows that $w_\alpha \to w_\infty$ with the rate of convergence
$\BIGO(\alpha^{\lambda_1/2})$ as $\alpha \to \infty$. This argument can be
made  rigorous although the details become non-trivial.

For $\alpha\neq \bar\alpha_n$, $\alpha>0$, we define an index
\[
  i(\alpha) = \# \left\{ r\in(0,\infty)\SEP w_\infty(r) = \bar\alpha_n\right\}\PERIOD
\]
The behaviour of $i(\alpha)$ in a neighbourhood of $\bar\alpha_n$ is controlled by the
linearization of the equation \VIZ{SE} at $w_{\bar\alpha_n}$. 
Denoting  $Z=\frac{\partial w_\alpha}{\partial \alpha}|_{\alpha=\bar\alpha_n}$ we have 
\begin{equation}\label{EZZ}
   Z'' + \left[\frac{n+1}{r} + (3\bar\alpha_n -1) r\right] Z' + 2Z = 0\COMMA
\end{equation}
with the initial conditions
\begin{equation}\label{ICZ}
  Z(0) = 1\COMMA\;\;\;Z'(0)=0\PERIOD
\end{equation}
The substitution $x = (1-3\bar\alpha_n)r^2/2$ transforms \VIZ{EZZ} into the 
standard form of the confluent hyper-geometric equation
\begin{equation}\label{CHG}
  x\frac{d^2 Z}{d x^2} + \left(\frac{n+2}{2} - x\right) Z + \frac{n+2}{n-4} Z = 0\PERIOD
\end{equation}
Hence the solution of \VIZ{EZZ} with the initial conditions \VIZ{ICZ} 
is found explicitly in the form
\[
  Z(r) = M\left(-\frac{n+2}{n-4}, \frac{n+2}{2}, \frac{n-4}{n+2}\frac{r^2}{2}\right)\COMMA
\]
where the function $M$ is one of the fundamental solutions of \VIZ{CHG} 
(see, for example, \cite{SpecFun}).
Using properties of the function $M$ and recalling that 
$\nu(n) = \min\{k\in\ZNUM\SEP (n+2)/(n-4)\leq k\}$ we see that
the solution $Z$ has $\nu(n)$ zeros in $(0,\infty)$. 
Before stating the following lemma we recall, that $R_\alpha$ denotes the first zero of
the solution $w_\alpha$, and that for $\alpha\neq\bar\alpha_n$, $\alpha>0$ 
$i(\alpha)=\#\{r,w_\alpha(r)=\bar\alpha_n\}$.

\begin{lemma}\label{lemma2}
  Assume $\bar\alpha_1<\alpha_1<\alpha_2$ or $0<\alpha_1<\alpha_2<\bar\alpha_n$,
  $R_{\alpha_1}$, $R_{\alpha_2} < \infty$, and $i(\alpha_1)\neq i(\alpha_2)$. Then there
  exists $\alpha\in(\alpha_1,\alpha_2)$ such that $R_\alpha=\infty$ and $w_\alpha$
  satisfies \VIZ{SE}, \VIZ{BC1}-\VIZ{BC2}.
\end{lemma}

\noindent
\PROOF To prove the lemma we use standard arguments based on the  continuity
of $i(\alpha)$ at points where $R_\alpha$ is finite. Some  work is required
to demonstrate that the condition \VIZ{BC2} is satisfied, but the arguments
 are straightforward.
 
\medskip 

From the behaviour of the linearized solution at $w_{\bar\alpha_n}$ one expects that
$R_\alpha$ is finite near $\bar\alpha_n$, $\alpha\neq\bar\alpha_n$ and
\begin{eqnarray} 
&&\lim_{\alpha\to\bar\alpha_n^+} = \left\{ \begin{array}{ll}
                                              \nu(n)     & \;\mbox{if $\nu(n)$ is odd} \\
                                              \nu(n) +1  & \;\mbox{if $\nu(n)$ is even} 
                                            \end{array}\right.   \label{IC1a}\\
&&\mbox{and} \\
&&\lim_{\alpha\to\bar\alpha_n^-} = \left\{ \begin{array}{ll}
                                             \nu(n) +1& \;\mbox{if $\nu(n)$ is odd} \\
                                             \nu(n)   & \;\mbox{if $\nu(n)$ is even} 
                                            \end{array}\right.  \label{IC1b}
\end{eqnarray}

The behaviour of $i(\alpha)$ for large $\alpha$ can be estimated from \VIZ{ASY}.
The function $M(a_1,b_1,x)$ for $a_1=\lambda_1/2$, $b_1=1+(\lambda_1+n-4+3\beta_n)/2$
has exactly one root in $(0,\infty)$ and decays exponentially to negative infinity as $x\to \infty$.
This behaviour suggests, that the solution $w_\alpha$, for large $\alpha$,
will be ``lagging behind'' $w_\infty$, and will fall back to $\bar\alpha_n$ for a finite, 
but large value of $r$.
After that  the term $r w'$ will be large and cause
the solution to reach $0$ in finite time. This reasoning leads us to expect
\begin{equation}\label{IC2}
  \lim_{\alpha\to\infty} i(\alpha) = 3\PERIOD
\end{equation}
Similar considerations suggest that
\begin{equation}\label{IC3}
  \lim_{\alpha\to 0^+} i(\alpha) = 2\PERIOD
\end{equation}
In this case we use  similar heuristic arguments: if $\alpha$ is very close to zero, the solution 
will reach the equilibrium 
$\bar\alpha_n$ at a large ``time'' $r$. The term $rw'(r)$ is then already significant.
After passing through $\bar\alpha_n$, the increase of $w$ will eventually be stopped
by the non-linear terms and the solution will start returning to $\bar\alpha$. 
After it passes through $\bar\alpha_n$ again, it will reach zero due to the
large term $rw'(r)$.
Assuming \VIZ{IC1a}-\VIZ{IC3}
we see, that for $\nu(n)$ odd, $i(\alpha)$ will change from $\nu(n)$ to  $3$ as $\alpha$ moves 
from a small right neighbourhood of $\bar\alpha_n$ to infinity and from
$\nu(n)+1$ to 2 as  $\alpha$ moves 
from a small left neighbourhood of $\bar\alpha_n$ to zero.

Assuming that $i(\alpha)$ changes by $2$ we obtain $\nu(n)-2$ solutions.
Repeating the same argument for $\nu(n)$ even yields again $\nu(n)-2$ of solutions.

The above described behaviour and the conjectured number 
of solutions are fully confirmed  by numerical computations. 
In Table~\ref{table1} we list
the solutions for integral values of $n$.
%
\begin{table}
\caption{Values of $\alpha$ for individual singular solutions tabulated
     at integral dimensions $n$. The branches of solutions, when continued in $n$,
         terminate at the dimensions indicated in the last column. 
         There are no solutions with $i(\alpha)=1$ or $i(\alpha)>6$.}\label{table1} 
\begin{indented}
\item[] \begin{tabular}{@{}lllllll}
\br
          & $n=5$ & $n=6$ & $n=7$ & $n=8$ & $n=9$ & max. $n$ \\
\mr
$i(\alpha)=2$&\NUM{0.02631}&\NUM{0.09647}&\NUM{0.1466}&\NUM{0.1684}&\NUM{0.17223}& $n=10$ \\
$i(\alpha)=3$&\NUM{1.3830}&\NUM{0.2792}&\NUM{0.2222}& &  & $n=7$  \\
$i(\alpha)=4$&\NUM{0.2205}&\NUM{0.2505}&     &       & &   $n=6$  \\
$i(\alpha)=5$&\NUM{0.2940}&  &       &       &       &     $n=11/2$ \\
$i(\alpha)=6$&\NUM{0.2855}&  &       &       &       &     $n=26/5$ \\
\br    
\end{tabular}
\end{indented}
\end{table}

\section{Semi-linear heat equation}\label{heq}
As we remarked in the introduction, there are striking similarities between the 
singular behaviour of solutions to \VIZ{eq4} when $n>4$ and the behaviour of
singular solutions of
\begin{equation}\label{HE}
   v_t = \LAP v + v^{2\sigma+1}\COMMA\;\;\;\mbox{in $\RR{n}\times(t_1,t_2)$,}
\end{equation}
for $n>10$ and $\sigma>\sigma_c(n)\equiv 2/(n-4-2\sqrt{n-1})$. For example, blow-up solutions
of \VIZ{HE} with the blow-up rate $(T-t)^{-1/\sigma-\delta}$ (with some $\delta>0$ 
for $n>10$ and  $\sigma>\sigma_c(n)$) were constructed in \cite{He-Va}. Such rate of blow-up
can be viewed as ``slower'' than the self-similar rate $(T-t)^{-1/\sigma}$. Some
authors call this rate ``faster'' but it seems that the term ``slower'' is more widely used
in the present context. It reflects the fact that $(T-t)^{-1/\sigma-\delta}$ becomes infinite
``more gradually''. 
We conjecture that the analysis of \cite{He-Va} can be used for proving
existence of slow blow-up solutions for $n>4$. In fact, if one allows slow decay of 
initial data at infinity, one can probably construct slow blow-up solutions also for $n\leq 4$.

In this section we look at radial self-similar singular solutions of \VIZ{HE}, i.e., 
the solutions of a special form
\begin{equation}\label{HSS}
   v(x,t) = \frac{1}{(2\kappa(T-t))^{1/\sigma}} w\left(\frac{|x|}{\sqrt{2\kappa(T-t)}}\right)\COMMA
\end{equation}
where $\kappa>0$, $T\in \RNUM$ are parameters. 
The notation of this section is not necessarily connected to the notation in 
the previous sections where we analyzed the equation \VIZ{eq1}.
Substituting \VIZ{HSS} into \VIZ{HE} we obtain the  equation for self-similar
profile
\begin{equation}\label{HSE}
   w'' +\frac{n-1}{r}w' - \kappa r w' + w^{2\sigma+1} -\frac{\kappa}{\sigma} w = 0\COMMA\;\;
\mbox{in $(0,\infty)$.}
\end{equation}
Natural boundary conditions are
\begin{eqnarray}
   && w(0)=\alpha\COMMA\;\;\; w'(0) = 0\COMMA\;\;\;\mbox{and} \label{HBC1}\\
   && w(r) \sim r^{-1/\sigma} \;\;\;\mbox{as $r\to\infty$.} \label{HBC2}
\end{eqnarray}
The equation \VIZ{HSE} has been studied by many authors and we recommend
the reader to consult \cite{Budd,Kohn-Giga,Lepin2,Lepin,Troy}.
Our aim in this section is to explain what happens to the solutions constructed
in \cite{Budd,Troy} and \cite{Lepin2,Lepin} as we approach the critical exponent
$\sigma_l(n) \equiv 3/(n-10)$ due to Lepin.

We use a similar approach as we applied to \VIZ{eq1} in the previous sections.
We shall denote $w_\alpha$ the solution of the initial-value problem \VIZ{HSE},
\VIZ{HBC1} for $\kappa=1$. We also define 
$R_\alpha = \sup\{0<r<\infty\SEP w_\alpha(r)>0\COMMA 0<r<\infty\}$.
Furthermore, we introduce the following parameters
\[
   \bar\alpha_\sigma =\left(\frac{1}{\sigma}\right)^{\frac{1}{2\sigma}}\COMMA
   \;\;\;\mbox{and}\;\;\;\;
   \beta_{n,\sigma} = \frac{1}{\sigma}\left(n-2-\frac{1}{\sigma}\right)\PERIOD
\]
As in the previous section it is convenient to use
new variables $(s,u)$ which are defined by the transformation
\[
   w(r) = r^{-1/\sigma}u(r)\COMMA\;\;\; r= e^s\PERIOD
\]
In these variables we obtain
\begin{equation}\label{HSEE}
   u'' + \left(n-2-\frac{2}{\sigma}\right) u' - \frac{1}{\sigma}\left(n-2-\frac{1}{\sigma}\right)u+
   u^{2\sigma+1} = \kappa e^{2s} u'\PERIOD
\end{equation}
The linearization of \VIZ{HSEE} at the trivial equilibrium $u=\beta_{n,\sigma}$ is 
\[
   z'' + \left( n-2-\frac{2}{\sigma}\right) z' + 
   2\left(n-2-\frac{1}{\sigma}\right) z = \kappa e^{2s}z'\PERIOD
\]
For the sake of brevity we denote $A= n-2-2/\sigma$ and $B=2(n-2-1/\sigma)$.
We denote $\lambda_1$, $\lambda_2$ the roots of the characteristic polynomial for 
the linear second-order differential operator on the left-hand side.
We note that the condition $A\geq 0$ is equivalent to $\sigma\leq\sigma_s(n)\equiv 2/(n-2)$.
It is known (see \cite{Kohn-Giga}) that \VIZ{HSE}, \VIZ{HBC1}-\VIZ{HBC2} has no bounded
solutions on $(0,\infty)$ in this case. We shall therefore assume $A<0$ 
(and hence $\mathrm{Re}\,\lambda_j<0$) in what follows. We also note that $\sigma\geq\sigma_c(n)$
corresponds to the requirement that $\lambda_j$ be real.

We define  $U(\alpha,\kappa,s)$ as the solution of the boundary-value problem given
by \VIZ{HSEE} and \VIZ{HBC1}-\VIZ{HBC2} rewritten in the variables $(s,u)$.
A formal calculation, similar to the one leading to \VIZ{ASY} in the previous section,
now gives
\begin{equation}\label{HASY}
\fl  U(e^{2s'}\alpha,\kappa,s) = \beta_n - 
   \REAL\left(\gamma e^{\lambda_1s'}e^{\lambda_1s}
   M\left(\frac{\lambda_1}{2},1+\lambda_1+\frac{A}{2},\frac{1}{2}\kappa e^{2s}\right)
                                                            \right)+ \dots \COMMA
\end{equation}
as $s'\to\infty$ for a suitable $\gamma\in\CNUM$, which is real if $\lambda_1$ is real. We emphasize
that we always assume $\REAL\lambda_1 < 0$.

The linearization $Z$ at $w_{\bar\alpha_\sigma}$ solves the boundary-value problem
\begin{eqnarray*}
  && Z'' + \frac{n-1}{r}Z' - r Z' + 2 Z = 0 \\
  && Z(0)=1\COMMA\;\;\; Z'(0)=0\PERIOD
\end{eqnarray*}
This equation has the solution $Z(r) = M(-1,n/2,r^2/2) \equiv 1 - r^2/n$.

We can introduce again the index $i(\alpha)\in \ZNUM$
\[
  i(\alpha) = \#\left\{ r\in(0,\infty)\SEP w_\alpha(r) = \bar\alpha_\sigma\right\}\COMMA
\]
and jumps of $i(\alpha)$ can be used to locate the desired solutions numerically.
Numerical computations show that $i(\alpha)$ jumps between $1$ and $3$.
Another natural index in this case (used in \cite{Budd, Lepin2,Lepin,Troy}) is
\[
  j(\alpha) = \#\left\{r\in(0,\infty)\SEP w_{\alpha}(r) = \frac{\beta_{n,\beta}}{r^{1/\sigma}}\right\}
                    = \#\left\{s\in\RNUM\SEP U(\alpha,1,s) = \beta_{n,\sigma}\right\}\PERIOD
\]
One has $j(\alpha)=2$ for $\alpha$ close to $\alpha_\sigma$. 
The approximation \VIZ{ASY} together with some heuristic arguments similar
to those used for justification of \VIZ{IC2} suggest that
\[
  \lim_{\alpha\to\infty} j(\alpha) = 
   \left\{ \begin{array}{ll}
          +\infty & \mbox{when $\sigma_s(n) < \sigma < \sigma_c(n)$,} \\
          \nu(n,\sigma) & \mbox{when $\sigma\geq\sigma_c(n)$ and $\nu(n,\sigma)$ is even,} \\
          \nu(n,\sigma)+1 & \mbox{when $\sigma\geq\sigma_c(n)$ and $\nu(n,\sigma)$ is odd,} 
                  \end{array}\right.
\]
where $\nu(n,\sigma)$ is defined as 
\[
  \nu(n,\sigma) = \min\left\{ k\in\ZNUM\SEP -\frac{\lambda_1}{2}\leq k\right\}\PERIOD
\]
This asymptotic behaviour of $j(\alpha)$ leads to infinitely many solutions for 
$\sigma_s(n) < \sigma < \sigma_c(n)$  (see \cite{Budd,Troy}).
When $R_{\alpha}<\infty$, $j(\alpha)$ must be even, hence typically we see $j(\alpha)$ 
jump by $2$. The solutions associated with these jumps have even index $j(\alpha)$.

However, the work of Lepin (\cite{Lepin}) shows that there are also solutions with odd 
index  $j(\alpha)$. Our numerical calculations suggest that such solutions can be detected
by jumps in the index $i(\alpha)$. The index $i(\alpha)$ takes on values $1$ or $3$.
The jump $1\to 3$ of the index $i(\alpha)$ when $\alpha$ increases corresponds to
the solutions with an even index $j(\alpha)$. On the other hand, the jump $3\to 1$ 
of $i(\alpha)$ indicates a solution with an odd index $j(\alpha)$. These solutions
cannot be numerically detected from the behaviour of $j(\alpha)$ only.

Lepin showed that solutions with the index $j(\alpha) = 2$ exist in the region
$\sigma > \sigma_s(n)$, $\nu(n,\sigma)\geq 3$. He also studied solutions
with higher indices, for which he established existence in smaller regions of
the parameter space $(n,\sigma)$. He conjectured that there are no
solutions to \VIZ{HSE}, \VIZ{HBC1}-\VIZ{HBC2} when $\nu(n,\sigma)\leq 2$,
which corresponds to $n> 10 + 3/\sigma$.

Our numerical computations strongly support this conjecture, at least 
in the sense that no solutions from the region $n < 10 + 3/\sigma$
can be continued outside the region. We looked at the solution with $j(\alpha)=2$
which appears to be the solution of \VIZ{HSE}, \VIZ{HBC1}-\VIZ{HBC2} with
the smallest possible $\alpha$. Numerically, it also is the most robust solution.
Assume $(\alpha_t,n_t,\sigma_t)$, $t\in [0,1)$ is a smooth path in the
parameter space such that $w_{\alpha_t}$ solves \VIZ{HSE}, \VIZ{HBC1}-\VIZ{HBC2}
for $n=n_t$, $\sigma=\sigma_t$ and $\kappa=1$. We consider the case $j(\alpha_t)=2$.
and follow the path of solutions with  $n_t - 10 - 3/\sigma_t <0$ such that
$n_t - 10 - 3/\sigma_t\to 0$, $n_t\to n_1$, $\sigma_t\to\sigma_1$ as $t\to 1$.
Numerically, we observe that $\alpha_t\to +\infty$. Moreover, for $n-10-3/\sigma\geq 0$
one expects from \VIZ{HASY} that $R_\alpha<+\infty$ for large $\alpha$, and
we have not detected any solutions for numerically accessible $\alpha$'s.
Therefore our conclusion is that the numerical computations provide strong
evidence in favour of Lepin's conjecture.

\ack
The research was supported in part by grants DMS--9877055, and DMS--0200326
from the National Science Foundation. P.P. acknowledges hospitality of the School
of Mathematics at the University of Minnesota. 
\section*{References}
\bibliographystyle{plain}
\bibliography{sveeq}

\end{document}